\documentclass[11pt,reqno]{article}
\usepackage{amsfonts, amsmath, amssymb, amscd, amsthm}
\usepackage{hyperref,verbatim}
\allowdisplaybreaks

\numberwithin{equation}{section} \textheight=22cm \textwidth=16cm
\newcommand{\divv}{\mathrm{div}}
\newcommand{\pd}{\frac{\partial}{\partial t}\Big|_{t=0}}

\newcommand{\grad}{\nabla}

\theoremstyle{plain}\newtheorem{thm}{Theorem}
\theoremstyle{plain}\newtheorem{lem}{Lemma}
\theoremstyle{definition}
\theoremstyle{definition}\newtheorem*{undefn}{Definition}
\theoremstyle{plain}\newtheorem{prop}{Proposition}
\theoremstyle{plain}
\theoremstyle{remark}\newtheorem{rem}{Remark}
\title{Two Kazdan-Warner type identities for the renormalized volume coefficients 
and the the Gauss-Bonnet curvatures of a Riemannian metric}

\author{Bin Guo \and Zheng-Chao Han \and Haizhong Li\thanks{Supported by grants of NSFC-10971110.}
}
\date{}
\begin{document}

\maketitle

\begin{abstract}
\noindent In this note, we prove two Kazdan-Warner type identities involving $v^{(2k)}$,
the renormalized volume coefficients of a Riemannian manifold $(M^n,g)$, and $G_{2r}$, the
so-called Gauss-Bonnet curvature, and a conformal Killing vector field on $(M^n,g)$. In
the case when the Riemannian manifold is locally conformally flat,
 $v^{(2k)}=(-2)^{-k}\sigma_k$, $G_{2r}(g)=\frac{4^r(n-r)!r!}{(n-2r)!}\sigma_r$ and our  results reduce to
 earlier ones established by Viaclovsky in
\cite{V00b} and  the second author in \cite{H06}.
  \end{abstract}



\vskip 20pt
\section{Introduction}

In \cite{V00b} and \cite{H06}, the  following result was proved

\medskip\noindent
{\bf Theorem A} (\cite{V00b}, \cite{H06}) {\it Let $(M,g)$ be a
  compact Riemannian manifold of dimension $n\geq 3$, $\sigma_k(g^{-1}\circ A_g)$ be the
$\sigma_k$ curvature of $g$, and $X$ be a conformal Killing vector
  field on $(M,g)$. When $k\geq 3$, we also assume that $(M,g)$ is locally conformally flat, then}
\begin{equation}\label{eqn:1.1}
\int_M\langle X,\nabla \sigma_k(g^{-1}\circ A_g)\rangle dv_g=0.
\end{equation}
Recall that  on an $n$-dimensional Riemannian manifold $(M,g)$, $n\geq 3$,
the full Riemannian curvature  tensor $Rm$ decomposes as
\begin{equation}\label{1.2}
Rm=W_g\oplus \left(A_g\odot\; g\right)
\end{equation}
where $W_g$ denotes the Weyl tensor of $g$,
\begin{equation}\label{1.3}
A_g=\frac{1}{n-2}({\rm Ric_g}-\frac{R_g}{2(n-1)}g)
\end{equation}
denotes the Schouten tensor, and $\odot$ is the Kulkani-Nomizu wedge product.  Under a conformal change of
metrics $g_w=e^{2w}g$,  where $w$ is a smooth function over the manifold, the Weyl curvature changes
pointwise as $W_{g_w}=e^{2w}W_g$. Thus, essential information of the Riemannian curvature tensor under a
conformal change of metrics  is reflected by the change of the Schouten tensor.  One often tries to study the
Schouten tensor through studying the elementary symmetric functions $\sigma_k(g^{-1}\circ A_g)$ (which we
later denote as $\sigma_k(g)$) of the eigenvalues  of the Schouten tensor, called the $\sigma_k$ curvatures
of $g$, and studying how they deform under conformal change of metrics.

The following question  is natural in relation to Theorem A:

\medskip\noindent
{\bf Question}. {\it Can we generalize Theorem A without the condition
  ``locally conformally flat"  for all $k\geq 1$?}

\medskip
In this note, we give an affirmative answer to the  above question.
Renormalized volume coefficients, $v^{(2k)}(g)$, of a Riemannian metric
$g$, were introduced in the physics literature in the late 1990's in
the context of AdS/CFT correspondence---see \cite{G00} for a
mathematical discussion, and were shown in \cite{GJ07} to be equal to
$\sigma_k(g^{-1}A_g)$, up to a scaling constant, when $(M,g)$ is locally conformally flat.
In fact, in the normalization we are going to adopt,
\begin{equation}\label{1.4}
 v^{(2)}(g)=-\frac{1}{2}\sigma_1(g),\qquad
 v^{(4)}(g)=\frac{1}{4}\sigma_2(g).
\end{equation}
 For $k=3$, Graham and Juhl (\cite{GJ07}, page 5) have aslo listed the following formula for $v^{(6)}(g)$:
\begin{equation}\label{1.5}
 v^{(6)}(g)=-\frac{1}{8}[\sigma_3(g)+\frac{1}{3(n-4)}(A_g)^{ij}(B_g)_{ij}],
\end{equation}
 where
\begin{equation}\label{ba}
 (B_g)_{ij}:=\frac{1}{n-3}\nabla^k\nabla^lW_{likj}+\frac{1}{n-2}R^{kl}W_{likj}
\end{equation}
 is the {\it Bach} tensor of the metric.
Just as
$\int_M \sigma_k(g^{-1}\circ A_g)\,dv_g$ is conformally invariant when $2k=n$
and $(M,g)$ is locally conformally flat, Graham  showed in
\cite{G00} that $\int_M v^{(2k)}(g) \,dv_g$ is also conformally invariant
on a general manifold when $2k=n$. Chang and Fang showed in \cite{CF08} that, for
$n\neq 2k$, the Euler-Lagrange equations for the functional $\int_M v^{(2k)}(g) \,dv_g$
under conformal variations subject to the constraint $Vol_g(M)=1$ satisfies
$v^{(2k)}(g)=$ const., which is a generalized characterization for the
curvatures $\sigma_k(g^{-1}\circ A_g)$ when $(M,g)$ is locally conformally flat, as given by Viaclovsky
\cite{V00}.

In this note, we will first show that the curvatures $v^{(2k)}(g)$ will play the role of
$\sigma_k(g^{-1}\circ A_g)$ in \eqref{eqn:1.1} for a general manifold. We note that
Graham \cite{G00} also gives an explicit expression of $v^{(8)}(g)$, but the explicit
expression of $v^{(2k)}(g)$ for general $k$ is not known because they are algebraically
complicated (see page 3 of \cite{G00}). Thus the study of the $v^{(2k)}(g)$ curvatures
involves significant challenges not shared by that of $\sigma_k(g)$:  firstly, for $k\geq
3$, $v^{(2k)}(g)$ depends on derivatives of curvature of $g$--- in fact, for $k\geq 3$,
$v^{(2k)}(g)$ depends on derivatives of curvatures of order up to $2k-4$; secondly, the
$v^{(2k)}(g)$ are defined via an indirect highly nonlinear inductive algorithm (see
\cite{G00}). Despite these difficulties, we can use some properties of these
$v^{(2k)}(g)$ curvatures to prove the following

\begin{thm}\label{thm:main}
Let $(M,g)$ be a compact Riemannian manifold of dimension $n\geq 3$, $X$ be a
conformal Killing vector field on $(M^n,g)$. For $k\geq 1$, we have
$$
\int_M\langle X,\nabla v^{(2k)}(g)\rangle dv_g=0.
$$
\end{thm}

 \medskip\noindent
\begin{rem} From (1.4), we know that Theorem \ref{thm:main} is equivalent to Theorem A when $k=1,2$, or when $(M^n,g)$ is
locally conformally flat for $k\geq 3$.
\end{rem}

The second result involves the Gauss-Bonnet curvatures $G_{2r}(2r\le n)$, introduced by
H. Weyl in 1939, which is defined by (also see \cite{L})
\begin{equation}
G_{2r}(g)=\delta^{j_1j_2\ldots j_{2r-1}j_{2r}}_{i_1i_2\ldots
i_{2r-1}i_{2r}}R^{i_1i_2}_{\quad j_1j_2}\ldots R^{{i_{2r-1}i_{2r}}}_{\qquad\;
j_{2r-1}j_{2r}},
\end{equation}
where $\delta^{j_1j_2\ldots j_{2r-1}j_{2r}}_{i_1i_2\ldots i_{2r-1}i_{2r}}$ is the generalized Kronecker
symbol. Note that $G_{2}=2R$, $R$ the scalar curvature. 
  We can prove that
\begin{thm}\label{thm:thm2}
Let $(M^n,g)$ be a compact Riemannian manifold, and $X$ be a conformal Killing vector
field. Then for the Gauss-Bonnet curvatures defined above, we have
\begin{equation}
\int _M \mathcal{L}_X G_{2r}(g)dv_g=0.
\end{equation}
\end{thm}
\begin{rem}
When $(M,g)$ is locally conformally flat, we see that the Gauss curvature
$G_{2r}(g)=\frac{4^r(n-r)!r!}{(n-2r)!} \sigma_r$, so Theorem \ref{thm:thm2} reduces to Theorem A.
\end{rem}
\begin{rem}
M. Labbi (\cite{L}) proved that the first variation of the functional $\int_M G_{2r}dv_g$
within the metrics with constant volume gave the so-called generalized Einstein metric,
and this functional has the variational property for $2r<n$ and is a topological
invariant for $2r=n$. In fact, if $n=2r$, this functional is the Gauss-Bonnet integrand
up to a constant (\cite{C}).
\end{rem}

 In the next section, we first provide a general proof for Theorem \ref{thm:main} by adapting an
ingredient in a preprint version of \cite{H06}, and making use of a variation formula for
$v^{(2k)}(g)$ established in \cite{G00} and \cite{CF08}.  And because of the explicit
expression for $v^{(6)}(g)$ and potential applications to other related problems in low
dimensions, we provide a self-contained proof for Theorem \ref{thm:main} in the case
$k=3$ in section 3. We will give a proof of Theorem \ref{thm:thm2} in section 4.

\section{Proof of Theorem \ref{thm:main}}\label{sec:2}

We will need the following
variation formula for $v^{(2k)}(g)$, see \cite{G00}.
\begin{prop}\label{prop:1}
Under the conformal transformation $g_t=e^{2t\eta}g$, the
variation of $v^{(2k)}(g_t)$ is given by
\begin{equation}\label{eqn:variation}
\pd v^{(2k)}(g_t)=-2k\eta v^{(2k)}+\nabla_i(L^{ij}_{(k)}\eta_j),
\end{equation}
where $L^{ij}_{(k)}$ is define as in \cite{G00} by
\begin{equation*}
L^{ij}_{(k)}=-\sum_{l=1}^k \frac{1}{l!}v^{(2k-2l)}(g) \partial^{l-1}_{\rho} g^{ij}(\rho) \Big |_{\rho=0},
\end{equation*}
with $g_{ij}(\rho)$ denoting the extension of $g$ such that
\[
g_{+}= \frac{(d\rho)^2-2\rho g(\rho)}{4\rho^2}
\]
is an asymptotic solution to $Ric(g_+)=-ng_+$ near $\rho=0$.
\end{prop}
An integral version of \eqref{eqn:variation} appeared  in \cite{CF08}:
\begin{equation}\label{3.2}
\int_M \big\{\frac{d}{dt}\Big |_{t=0}[v^{(2k)}(g_t)]+2k\eta v^{(2k)}(g)\big\}dv_g=0.
\end{equation}
\begin{proof}[Proof of Theorem~\ref{thm:main} in the case $n\neq 2k$.]
 Let $X$ be a conformal vector field on $M$. Let $\phi_t$ denote the local one-parameter family of
 conformal diffeomorphisms of $(M,g)$ generated by $X$. Thus for some smooth function $\omega_t$ on $M$, we have
 \begin{equation}\label{eqn:26}
 \phi^*_t(g)=e^{2\omega_t}g=:g_t.
 \end{equation}
 We have the following properties
\begin{equation}\label{eqn:v}
\phi_t^* v^{(2k)}(g)=v^{(2k)}(\phi_t^*g)=v^{(2k)}(e^{2\omega_t}g),
\end{equation}
\begin{equation}\label{eqn:28}
 {\dot\omega}:=\frac{d}{dt}\Big |_{t=0}\omega_t=\frac{{\rm div}X}{n},
\end{equation}
 \begin{equation}\label{eqn:29}
 \frac{d}{dt}\Big |_{t=0}\left(g_t^{-1}\circ A(g_t)\right)=-\grad^2{\dot\omega}-2{\dot\omega}g^{-1}\circ A(g).
 \end{equation}
\begin{equation}\label{vdiv}
\frac{\partial}{\partial t}\Big|_{t=0}\divv_{g_t}X=n
X\eta=n\langle X,\nabla \eta\rangle.
\end{equation}
Using \eqref{eqn:v}, \eqref{eqn:28}, and \eqref{eqn:variation},
we have
\begin{equation*}\begin{split}
\langle X,\nabla v^{(2k)}(g)\rangle &= \frac{d}{dt}\Big |_{t=0}[v^{(2k)}(g_t)] \\
&=-2k\dot{\omega} v^{(2k)}+\nabla_i(L^{ij}_{(k)}\dot{\omega}_j)\\
&=-\frac{2k}{n} (\divv X)  v^{(2k)}+\nabla_i(L^{ij}_{(k)}\dot{\omega}_j)\\
&= - \frac{2k}{n}\divv ( v^{(2k)} X)+ \frac{2k}{n} \langle X,\nabla v^{(2k)}(g)\rangle
+\frac 1n \nabla_i(L^{ij}_{(k)}(\divv X)_j),
\end{split}\end{equation*}
from which it follows that
\begin{equation}\label{han06}
\left(1-  \frac{2k}{n}\right) \langle X,\nabla v^{(2k)}(g) \rangle
= - \frac{2k}{n}\divv ( v^{(2k)} X) + \frac 1n \nabla_i(L^{ij}_{(k)}(\divv X)_j).
\end{equation}
Theorem \ref{thm:main} in the case $2k \ne n$ now follows directly by integrating \eqref{han06} over $M$.
\end{proof}
\begin{proof}[Proof of Theorem \ref{thm:main} in the case $2k=n$.]
As in \cite{H06}, we will prove that for any conformal metric $g_1=e^{2\eta}g$ of $g$,
\begin{equation}\label{eqn:be}
\int_M \langle X, v^{(2k)}(g_1)\rangle d v_{g_1} =\int_M \langle X, v^{(2k)}(g)\rangle d v_{g} =-\int_M {\rm
div}_g X v^{(2k)}(g) d v_g,
\end{equation}
 i.e. $\int_M\langle X, v^{(2k)}(g)\rangle dv_{g}$ is
independent of the particular choice of metrics in the conformal
class. To this end, we only have to prove that for
$g_t=e^{2t\eta}g$,
\begin{equation}\label{eqn:int}
\pd \int_M \divv_{g_t}X v^{(2k)}(g_t) dv_{g_t}=0.
\end{equation}
We prove \eqref{eqn:int} by direct computations using Proposition
\ref{prop:1}. Indeed,
\begin{equation}\label{eqn:v2k}\begin{split}
&\pd  \int_M \divv_{g_t}X v^{(2k)}(g_t) dv_{g_t}\\
=&\int_M \Big[n\langle X,\nabla\eta\rangle v^{(2k)}+\divv X\big(-2k\eta
v^{(2k)}+\nabla_i(L^{ij}_{(k)}\eta_j)\big)
+n\eta \divv X v^{(2k)}\Big] dv_g\\
=&\int_M \Big[n\langle X,\nabla\eta\rangle v^{(2k)}+\divv X \nabla_i(L^{ij}_{(k)}\eta_j)\Big]dv_g\\
=& \int_M \Big[\langle n v^{(2k)} X, \nabla\eta\rangle -L^{ij}_{(k)} (\divv X)_i \eta_j)\Big]dv_g\\
=& \int_M \Big[ - \divv ( n v^{(2k)} X)+\nabla_j\left(L^{ij}_{(k)} (\divv X)_i\right)\Big] \eta dv_g =0
\end{split}\end{equation}
in the case $n=2k$ by \eqref{han06}.

The remaining argument is an adaptation of
an argument of Bourguignon and Ezin ([BE]):
either the connected component of the identity of the conformal
group $C_0(M,g)$ is compact, then there is a metric $\hat{g}$
conformal to $g$ admitting $C_0(M,g)$ as a group of isometries, from
which it follows that $\divv _{\hat g}X\equiv 0$ and (1.7) therefore
holds; or, $C_0(M,g)$ is non-compact, then by a theorem of
Obata-Ferrand, $(M,g)$ is conformal to the standard sphere, in which
case we can pick the canonical metric to compute the integral on the
left hand side of $(1.7)$ and conclude that it is zero.
\end{proof}

\section{Self-contained proof of Theorem \ref{thm:main} in the case  $k=3$}

We aim to give a direct, self-contained derivation for a more explicit version of \eqref{eqn:variation}, more precisely,
under conformal change of metrics $g_{t}=e^{2t\eta}g$, we have
\begin{equation}\label{key6}
\pd v^{(6)}(g_t) =-6 v^{(6)}(g) \eta +\nabla^j
\bigg[\Big(\frac{T^{(2)}_{ij}(g)}{8}+\frac{B_{ij}(g)}{24(n-4)}\Big)\nabla^i \eta\bigg],
\end{equation}
where $T^{(2)}_{ij}(g)$ is the Newton tensor associated with $A_g$, as defined in Reilly [R]:
\begin{undefn}For an integer $k\geq 0$, $k$-th Newton tensor is
$$
T^{(k)}_{ij}=\frac{1}{k!}\sum\delta^{j_1\cdots j_k j}_{i_1\cdots
i_k i}A_{i_1j_1}\cdots A_{i_kj_k}
$$
where $\delta^{j_1\cdots j_k j}_{i_1\cdots
i_k i}$ is the generalized Kronecker symbol.
\end{undefn}
With \eqref{key6} we can repeat the proof in the last section to prove Theorem \ref{thm:main} in the case
$k=3$.

First we recall the transformation laws for the tensors $B_{ij}$ and
$A_{ij}$ under conformal change of metrics $g_{t}=e^{2t\eta}g$---see \cite{CF08}:
\begin{equation*}
A_{ij}(g_t)=A_{ij}-\nabla^2_{ij}\eta+\nabla_{i}\eta\nabla_{j}\eta-\frac{|\nabla \eta|_g^2}{2}g_{ij};
\end{equation*}
\begin{equation*}
B_{ij}(g_t)=e^{-2t\eta}\Big(B_{ij}+(n-4)t(C_{ijk}+C_{jik})\grad^k
\eta+(n-4)t^2W_{ikjl}\grad^k\eta \grad^l\eta\Big).
\end{equation*}
where $C_{ijk}$ are the components of the  {\it Cotton} tensor  defined  by
$$
C_{ijk}=A_{ij,k}-A_{ik,j}
$$
with $A_{ij,k}$ being the components of the covariant derivative of the Schouten tensor $A_{ij}$.

Thus
\begin{equation*}
\pd A^{ij}(g_t) = - \grad^{ij}\eta-4A^{ij}(g)\eta,
\quad \text{and} \quad \pd B_{ij}(g_t) = (n-4)(C_{ijk}+C_{jik})\grad^k \eta - 2\eta B_{ij}.
\end{equation*}
We recall some properties to be used.
\begin{prop}\label{prop:prt} (\cite{V00},\cite{H06},[HL]). We have
\begin{enumerate}
\item[(i)] $k\sigma_{k}(g)=\sum\limits_{i,j}T^{(k-1)}_{ij}A_{ij}$
\item[(ii)] $\sum\limits_i T^{(k)}_{ii}=(n-k)\sigma_k(g)$.
\item[(iii)] $\sum\limits_l\nabla^l W_{lijk}=-(n-3)C_{ijk}$.
\end{enumerate}
\end{prop}
Using the relation between $v^{(6)}$ and $\sigma_3(g)$, $A^{ij}B_{ij}$ as in \eqref{1.5}, we find
\begin{equation*}
\begin{split}
&-8\pd v^{(6)}(g_t) \\
=& T^{(2)}_{ij}(g)\bigg(-\grad^{ij}\eta-2\eta A^{ij}(g)\bigg)+ \frac{1}{3(n-4)}\bigg[-B_{ij}(g)
\grad^{ij}\eta + (n-4) A^{ij}(g) (C_{ijk}+C_{jik})\grad^k\eta-
6 \eta A^{ij}B_{ij}\bigg]\\
=& -6 \bigg(\sigma_3(g)+ \frac{1}{3(n-4)} A^{ij}B_{ij}\bigg)\eta -\bigg[
T^{(2)}_{ij}(g)+\frac{B_{ij}(g)}{3(n-4)}\bigg]
\grad^{ij}\eta +\frac{2}{3}  A^{ij}(g) C_{ijk} \grad^k\eta\\
=& 48  v^{(6)}(g) \eta -\nabla^j \bigg[\Big(T^{(2)}_{ij}(g)+\frac{B_{ij}(g)}{3(n-4)}\Big)\nabla^i \eta\bigg]
+\bigg[\sum_j \Big(T^{(2)}_{ij,j}(g)+ \frac{B_{ij,j}(g)}{3(n-4)}\Big)+\frac{2}{3} A^{kl}C_{kli}\bigg]\grad^i
\eta,
\end{split}
\end{equation*}
where we used \eqref{1.5} and (i) of Proposition~\ref{prop:prt}.
In the following we will verify that
\[
\sum_j \left(T^{(2)}_{ij,j}(g)+ \frac{B_{ij,j}(g)}{3(n-4)}\right)+\frac{2}{3}  A^{kl}C_{kli}=0,
\]
thus establishing \eqref{key6}.
The above property would follow from the following
\begin{lem}\label{3.1}
\begin{enumerate}
\item[(i)] $\sum\limits_jT^{(2)}_{ij,j}=-A^{pq}C_{pqi}$;
\item[(ii)] $\sum\limits_j B_{ij,j}=(n-4)A^{kl}C_{kli}.$
\end{enumerate}
\end{lem}

\begin{proof}[Proof of (i).] We have the following calculation in normal coordinate,
\begin{equation*}
\begin{split}
\sum\limits_jT^{(2)}_{ij,j}&=\sum\limits(
\frac{1}{2!}\sum\delta^{j_1j_2 j}_{i_1i_2 i}A_{i_1j_1}A_{i_2j_2})_j\\
&=
\sum\delta^{j_1j_2 j}_{i_1i_2 i}A_{i_1j_1}A_{i_2j_2,j}\\
&= -A^{pq}C_{pqi},
\end{split}
\end{equation*} where we used
$$
\delta^{j_1j_2j}_{i_1i_2i}= \left |
\begin{array}{cccc}
\delta_{i_1j_1}&\delta_{i_1j_2}&\delta_{i_1j}\\
\delta_{i_2j_1}&\delta_{i_2j_2}&\delta_{i_2j}\\
\delta_{ij_1}&\delta_{ij_2}&\delta_{ij}
\end{array}
\right | 
$$
and $\sum\limits_iA_{ii,j}=\sum\limits_i A_{ij,i}$, which itself is a consequence of the second
Bianchi identity.
\end{proof}
\medskip\noindent
\begin{proof}[Proof of (ii).]
First, using (iii) of Proposition \ref{prop:prt} and substituting $R_{ij}$ in terms of $A_{ij}$ in the
definition of the Bach tensor $B_{ij}$, we obtain
\begin{equation*}
\begin{split}
B_{ij}&=-\sum_k C_{ikj,k}+\sum_{k,l}A_{kl} W_{likj}\\
&= -\sum_k \left( A_{ik,jk }-A_{ij,kk}\right)+ \sum_{k,l}A_{kl} W_{likj}.
\end{split}
\end{equation*}
Thus
\begin{align*}
&\sum_j B_{ij,j}\\
=& - \sum_{j,k} \left( A_{ik,jkj}-A_{ij,kkj} \right)+ \sum_{k,l,j}
\left(A_{kl,j} W_{likj}+ A_{kl}W_{likj,j}\right)\\
=& - \sum_{j,k} \left( A_{ik,jkj}-A_{ik,jjk} \right)+ \sum_{k,l,j}A_{kl,j} W_{likj}-(n-3)\sum_{k,l}A_{kl}C_{kil}\\
=&- \sum_{j,k,m} \left(A_{ik,m}R_{mjkj}+A_{im,j}R_{mkkj}+A_{mk,j}R_{mikj}\right)+
 \sum_{k,l,j}A_{kl,j} W_{likj}+(n-3)\sum_{k,l}A_{kl}C_{kli}\\
=& \sum_{j,k,m} \left(-A_{mk,j}R_{mikj} + A_{km,j}W_{mikj}\right)+(n-3)\sum_{k,l}A_{kl}C_{kil}\\
=& \sum_{j,k,m}A_{mk,j}\left(-A_{mk}g_{ij}+A_{mj}g_{ik}-g_{mk}A_{ij}+g_{mj}A_{ik}\right)+(n-3)\sum_{k,l}A_{kl}C_{kli}\\
=& \sum_{m,k} \left(-A_{mk,i} A_{mk}+A_{mi,k}A_{mk}-A_{mk,j}g_{mk}A_{ij}+A_{mj,k}g_{mk}A_{ij}\right)
+(n-3)\sum_{k,l}A_{kl}C_{kli}\\
=& \sum_{m,k}  A_{mk}\left(A_{mi,k}-A_{mk,i}\right)
+(n-3)\sum_{k,l}A_{kl}C_{kli}\\
=& \sum_{m,k} A_{mk} C_{mik}
+(n-3)\sum_{k,l}A_{kl}C_{kli}\\
=& (n-4) \sum_{k,l}A_{kl}C_{kli},
\end{align*}
where we have used  
\[
R_{mikj}=W_{mikj}+ A_{mk}g_{ij}-A_{mj}g_{ik}+g_{mk}A_{ij}-g_{mj}A_{ik}.
\]
\end{proof}
\begin{proof}[Proof of Theorem \ref{thm:main} of the special case $k=3$]
We use the notations of section 2, let $\phi_t$ be the local one-parameter family of conformal
diffeomorphisms of $(M,g)$ generated by $X$. For $g_t=\phi_t^*(g)=e^{2\omega_t}g$, similar to \eqref{key6} we
have
\begin{equation}\label{eqn:main}
\langle X,v^{(6)}\rangle=\pd v^{(6)}(g_t) =-6 v^{(6)}(g) \dot\omega +\sum\limits_{i,j}\nabla^j
\bigg[\Big(\frac{T^{(2)}_{ij}(g)}{8}+\frac{B_{ij}(g)}{24(n-4)}\Big)\nabla^i \dot\omega\bigg],
\end{equation}
if $n\neq 2k$ then integrating \eqref{eqn:main} we can get Theorem \ref{thm:main}.

While if $n=2k$, then by use of \eqref{key6} and \eqref{eqn:main}, we can prove that $\int_M\langle X,
v^{(6)}(g)\rangle dv_{g}$ is independent of the particular choice of the metric within the conformal class.
The remaining of the proof is verbatim the same as that of section 2.
\end{proof}
\section{Proof of Theorem \ref{thm:thm2}}
In this section, we will prove Theorem \ref{thm:thm2} using a similar method as in section
\ref{sec:2}. Let $(M^n,g)$ be a compact Riemannian manifold, and we denote by $R_{ijkl}$
the Riemann curvature tensor in local coordinates. Define a tensor $P_r (2r\le n)$ by
\begin{equation*}
{P_r}_{i}^{\; j}= \delta^{jj_1j_2\ldots j_{2r-1}j_{2r}}_{ii_1i_2\ldots
i_{2r-1}i_{2r}}R^{i_1i_2}_{\quad j_1j_2}\ldots R^{{i_{2r-1}i_{2r}}}_{\qquad\;
j_{2r-1}j_{2r}},
\end{equation*}
where $\delta^{jj_1j_2\ldots j_{2r-1}j_{2r}}_{ii_1i_2\ldots i_{2r-1}i_{2r}}$ is the
generalized Kronecker symbol. First we give the following lemma.
\begin{lem}\label{lem:lem1} The tensor $P_{r}$ is divergence free, i.e.
$${P_r}_{i,j}^j=0, \text{ for any $i$.}$$
\end{lem}
\begin{proof}
We have the following direct computations.
\begin{equation*}\begin{split}
{P_r}_{i,j}^j&=r\delta^{jj_1j_2\ldots j_{2r-1}j_{2r}}_{ii_1i_2\ldots
i_{2r-1}i_{2r}}R^{i_1i_2}_{\quad
j_1j_2,j}\ldots R^{{i_{2r-1}i_{2r}}}_{\qquad\; j_{2r-1}j_{2r}}\\
&=-r\delta^{jj_1j_2\ldots j_{2r-1}j_{2r}}_{ii_1i_2\ldots i_{2r-1}i_{2r}}R^{i_1i_2}_{\quad
j_2j,j_1}\ldots
R^{{i_{2r-1}i_{2r}}}_{\qquad\; j_{2r-1}j_{2r}}\\
& \quad - r\delta^{jj_1j_2\ldots j_{2r-1}j_{2r}}_{ii_1i_2\ldots
i_{2r-1}i_{2r}}R^{i_1i_2}_{\quad jj_1,j_2}\ldots
R^{{i_{2r-1}i_{2r}}}_{\qquad\; j_{2r-1}j_{2r}}\\
&= -2r\delta^{jj_1j_2\ldots j_{2r-1}j_{2r}}_{ii_1i_2\ldots
i_{2r-1}i_{2r}}R^{i_1i_2}_{\quad j_1j_2,j}\ldots
R^{{i_{2r-1}i_{2r}}}_{\qquad\; j_{2r-1}j_{2r}}\\
&= -2{P_r}_{i,j}^j,
\end{split}\end{equation*}
where we have used the second Bianchi identity. It then follows that ${P_r}_{i,j}^j=0$.
\end{proof}
We need the following algebraic lemma.
\begin{lem}\label{lem:lem2}
The generalized Kronecker symbol satisfies
\begin{equation*}
\sum\limits_{i,j=1}^n\delta^i_j\delta_{ii_1\ldots i_r}^{jj_1\ldots j_r}=
(n-r)\delta^{j_1\ldots j_r}_{i_1\ldots i_r},
\end{equation*}
for any $1\le i_1,\ldots, j_r\le n$, and $r\le n$.
\end{lem}
The proof of Lemma \ref{lem:lem2} is a direct calculation by use of the definition and we
omit it here.

 Let $X$ be a conformal vector field, denoted by $\phi_t$ be the
one-parameter subgroup of diffeomorphism generated by $X$. Then there exists a family of
functions $\omega_t$ such that $g_t=\phi_t^* g=e^{2\omega_t}g$. We have \eqref{eqn:28},
$\omega_0=0$, and
\begin{equation}\label{eqn:eqn}
G_{2r}(g_t)=\phi_t^* G_{2r}(g).
\end{equation}
Under conformal change of metrics $g_t=e^{2\omega_t}g$, we have the following formula
(see e.g. \cite{CLN}),
\begin{equation}\label{eqn:rm}
R^{ij}_{\;\; kl}(g_t)=e^{-2\omega_t}\Big(R^{ij}_{\;\; kl}-(\alpha \odot g)^{ij}_{\;\;
kl}\Big),
\end{equation}
where we denote $\alpha_{ij}=({\omega_t})_{ij}-({\omega_t})_i({\omega_t})_j+\frac{|\nabla
\omega_t|^2}{2}g_{ij}$ for convenience (note that $({\omega_t})_{ij}$ is the covariant
derivative with respect to the fixed metric $g$.) and $\odot$ is the Kulkani-Nomizu
product, defined by
\begin{equation*}
(\alpha\odot
g)_{ijkl}=\alpha_{ik}g_{jl}+\alpha_{jl}g_{ik}-\alpha_{il}g_{jk}-\alpha_{jk}g_{il}.
\end{equation*}
 From \eqref{eqn:rm} we see
that
\begin{equation}\label{eqn:f2r}
G_{2r}(g_t)=e^{-2r\omega_t}\delta^{j_1j_2\ldots j_{2r-1}j_{2r}}_{i_1i_2\ldots
i_{2r-1}i_{2r}}\Big(R^{i_1i_2}_{\quad j_1j_2}-(\alpha\odot g)^{i_1i_2}_{\quad
j_1j_2}\Big)\ldots \Big(R^{i_{2r-1}i_{2r}}_{\qquad\quad j_{2r-1}j_{2r}}- (\alpha\odot
g)^{i_{2r-1}i_{2r}}_{\qquad\quad j_{2r-1}j_{2r}}\Big).
\end{equation}
Taking derivative with respect to $t$ on both sides of \eqref{eqn:eqn} and using
\eqref{eqn:f2r}, we see by use of \eqref{eqn:28}
\begin{equation}\label{eqn:evo}\begin{split}
\mathcal{L}_X G_{2r}(g)&=\frac{\partial}{\partial t}\Big|_{t=0}G_{2r}(g_t)\\
  &= -2r\dot\omega G_{2r}(g)- r \delta^{j_1j_2\ldots j_{2r-1}j_{2r}}_{i_1i_2\ldots i_{2r-1}i_{2r}}
  \big(\frac{\partial\alpha}{\partial t}\Big|_{t=0}\odot g\big)^{i_1i_2}_{\quad j_1j_2}R^{i_3i_4}_{\quad j_3j_4} \ldots
  R^{i_{2r-1}i_{2r}}_{\qquad\quad j_{2r-1}j_{2r}}\\
  &= -2r\dot\omega G_{2r}(g) - 4r(n-2r+1) {P_{r-1}}^{\; j}_{i}\dot\omega^i_j\\
  &= -2r\frac{\mathrm{div}X}{n}G_{2r}(g)-\frac{4r(n-2r+1)}{n}{P_{r-1}}^{j}_i (\mathrm{div}X)^i_{\;j}\\
  &= -2r\frac{\mathrm{div}X}{n}G_{2r}(g)-\frac{4r(n-2r+1)}{n}\nabla_j\Big({P_{r-1}}^{j}_i
  (\mathrm{div}X)^i\Big).
\end{split}\end{equation}
where we have used Lemma \ref{lem:lem2} in the third equality and Lemma \ref{lem:lem1} in
the last equality. Integrating \eqref{eqn:evo} over $M$ and using the divergence theorem,
we see that
\begin{equation}\label{eqn:eqmain}\begin{split}
\int_M \mathcal{L}_X G_{2r}(g) dv& =-2r \int_M \frac{\mathrm{div}X}{n}G_{2r}(g) dv
=\frac{2r}{n} \int_M \mathcal{L}_X G_{2r}(g)dv,
\end{split}\end{equation}
Hence, if $n>2r$, it follows from \eqref{eqn:eqmain} that $\int_M \mathcal{L}_X G_{2r}(g) dv
= 0$. If $n=2r$, we follow similar ideas as in section \ref{sec:2}, i.e. we need to prove
that the integral
\begin{equation*}
\int _M G_{2r}(g) \mathrm{div}_g X dv_g,
\end{equation*}
is independent of a particular choice of  metrics within a conformal class. Let
$g_1=e^{2\eta}g(\eta\in C^\infty(M))$ be any metric in the conformal class $[g]$.
Considering a family of metrics $g_t=e^{2t\eta}g $ connecting $g$ and $g_1$, we need to
prove that
\begin{equation*}
\frac{d}{dt}\Big|_{t=0}\int _M G_{2r}(g_t) \mathrm{div}_{g_t} X dv_{g_t}=0.
\end{equation*}
By a direct computation, we have
\begin{equation*}\begin{split}
&\frac{d}{dt}\Big|_{t=0}\int _M G_{2r}(g_t) \mathrm{div}_{g_t} X dv_{g_t}\\
=&\int_M\Big[ \frac{\partial}{\partial t}\Big|_{t=0}G_{2r}(g_t)\mathrm{div}X+
G_{2r}(g)\frac{\partial}{\partial
t}\Big|_{t=0} \mathrm{div}_{g_t}X + n \eta G_{2r}(g)\mathrm{div}X\Big]dv_g\\
=& \int_M \Big[-2r\eta G_{2r}(g)\mathrm{div}X -4r(n-2r+1)
{P_{r-1}}^{\;j}_i\eta^i_{\;j}\mathrm{div}X +
nG_{2r}(g)\langle\nabla\eta, X\rangle + nG_{2r}(g)\mathrm{div}X \eta\Big]dv_g\\
=& \int_M \Big[-2r\eta G_{2r}(g)\mathrm{div}X -4 \eta r(n-2r+1){P_{r-1}}^{\;j}_i
(\mathrm{div}X)^i_j- n\eta
\langle \nabla G_{2r}(g), X \rangle\Big]dv_g\\
=& 0,
\end{split}\end{equation*}
where we have used \eqref{vdiv} in the second equality, the divergence theorem in the
third equality and \eqref{eqn:evo} in the last equality. The remaining proof follows the
idea of \cite{BE} as in section \ref{sec:2}. Hence we complete the proof of Theorem
\ref{thm:thm2}.

\vskip 1cm
\begin{flushleft}
\medskip\noindent

Bin Guo: {\sc Department of Mathematical Sciences, Tsinghua
University, Beijing 100084, People's Republic of China}\ \  Email:
guob07@mails.tsinghua.edu.cn

Zheng-Chao Han: {\sc Department of Mathematics, Rutgers
University, 110, Frelinghuysen Road, Piscataway, NJ 08854, USA} \ \ E-mail:
zchan@math.rutgers.edu

Haizhong Li: {\sc Department of Mathematical Sciences, Tsinghua
University, Beijing 100084, People's Republic of China} \ \ E-mail:
hli@math.tsinghua.edu.cn

\end{flushleft}
\end{document}